\documentclass{article}
\usepackage[english]{babel}
\usepackage{amsthm}
\usepackage[dvipsnames]{xcolor}
\usepackage{amsmath}
\usepackage{amssymb}
\usepackage[%
]{hyperref}
\usepackage{booktabs}
\usepackage{enumerate}
\usepackage[left = 20mm, top = 35mm, bottom = 30mm, right = 20mm]{geometry}
\usepackage{multicol}
\usepackage{bbm}
\usepackage[]{todonotes}
\usepackage{setspace}
\usepackage{mathtools}
\usepackage{tikz}
\usepackage{subcaption}

\tikzset
{
	edge/.style = {-, very thick},
	every node/.style = {circle, draw, fill, inner sep=2pt, minimum size=6pt},
}

\theoremstyle{definition}
\newtheorem{definition}{Definition}
\newtheorem{Remark}[definition]{Remark}
\newtheorem{Example}[definition]{Example}

\theoremstyle{plain}
\newtheorem{theorem}[definition]{Theorem} 
\newtheorem{theoremintro}{Theorem}
\newtheorem{lemma}[definition]{Lemma}

\newcommand{\N}{\mathbb{N}}

\newcommand{\mB}{\mathcal{B}}
\newcommand{\mR}{\mathcal{R}}
\newcommand{\rC}{M_R}
\newcommand{\bC}{M_B}
\newcommand{\Aut}{\operatorname{Aut}}

\newcommand{\myred}{Melon}
\newcommand{\myblue}{Cyan}

\newcommand\extrafootertext[1]{%
	\bgroup
	\renewcommand\thefootnote{\fnsymbol{footnote}}%
	\renewcommand\thempfootnote{\fnsymbol{mpfootnote}}%
	\footnotetext[0]{#1}%
	\egroup
}

\numberwithin{definition}{section}
\numberwithin{equation}{section}

\title{Countable ultrahomogeneous graphs on two imprimitive color classes}
\author{\textbf{Sofia Brenner and Irene Heinrich} \\
	\normalsize{TU Darmstadt, Fachbereich Mathematik}\\
\normalsize{Schlossgartenstr.~7, 64289 Darmstadt, Germany}\\
\normalsize{\texttt{\{brenner, heinrich\}@mathematik.tu-darmstadt.de}}}
\date{\vspace{-0.8cm}}

\begin{document}
	\maketitle
	\begin{abstract}
		\noindent
		We classify the countable ultrahomogeneous $2$-vertex-colored graphs in which the color classes are imprimitive, i.e., up to complementation they form disjoint unions of cliques. This generalizes work by Jenkinson et.~al.~\cite{JEN12}, Lockett and Truss \cite{LOC14} as well as Rose~\cite{ROS11} on ultrahomogeneous $n$-graphs. As the key aspect in such a classification, we identify a concept called piecewise ultrahomogeneity. We prove that there are two specific graphs whose occurrence essentially dictates whether a graph is piecewise ultrahomogeneous, and we exploit this fact to prove the classification.
		\bigskip
		
		
	\end{abstract}

	\section{Introduction}
	\extrafootertext{An extended abstract of this paper~\cite{BRE23} is to appear in the Proceedings of the European Conference on Combinatorics, Graph Theory
		and Applications (EUROCOMB’23). This research has received funding from the European Research Council (ERC) under the European Unions Horizon 2020 research and innovation programme (EngageS: grant agreement No.~820148). We thank Pascal Schweitzer for helpful discussions.}
	Ultrahomogeneous structures are relational structures in which every isomorphism between finite substructures can be extended to an automorphism of the entire structure.\footnote{Some authors use the term ``homogeneous'' for this property.} The extensive study of ultrahomogeneous objects relates various areas of research, such as model theory, permutation group theory and Ramsey theory (see \cite{MAC11} for a survey). A vast collection of ultrahomogeneous classes of relational structures has been classified. For instance, apart from different classes of graphs, which we discuss below, there exist classification results for partially ordered sets \cite{SCH79}, tournaments \cite{LAC84, CHE98} as well as countably infinite permutations \cite{CAM02}.
	\medskip
	
	In this article, we focus on a special class of countable ultrahomogeneous graphs. By work of Sheehan \cite{SHE74} and Gardiner \cite{GAR76} as well as Golfand and Klin \cite{GOL78}, the finite ultrahomogeneous graphs are known. Lachlan and Woodrow~\cite{LAC80} gave a characterization of the ultrahomogeneous graphs with countably infinitely many vertices. Cherlin \cite{CHE98} asked for a classification of ultrahomogeneous $n$-graphs, that is, ultrahomogeneous graphs for which the vertex set is partitioned into $n$ subsets which are respected by the partial isomorphisms considered. 
	\medskip
	
	Nowadays, one usually thinks of $n$-graphs as graphs with a vertex-coloring in $n$ colors, and considers isomorphisms preserving colors. Finite ultrahomogeneous vertex-colored graphs were classified in \cite{HEI20}. Every color class in an ultrahomogeneous graph induces a monochromatic ultrahomogeneous graph. In particular, up to complementation every infinite color class forms an imprimitive graph (a nontrivial disjoint union of cliques), an independent set, or it induces a Rado graph or a Henson graph (see \cite{LAC80}). Jenkinson et.~al.~\cite{JEN12} considered vertex-colored graphs in which the color classes form independent sets. Their work was extended by Lockett and Truss \cite{LOC14} who allowed an additional coloring of the edges (while still requiring that every color class forms an independent set). In his dissertation, Rose \cite{ROS11} investigates countable 2-colored graphs. The main part of his work covers the case that one color class forms a disjoint union of cliques and the other one induces a Rado graph or a Henson graph. For the case that both color classes form a disjoint union of cliques, a partial list of possible cases is stated, but not proven. 
	\medskip

	In this paper, we classify the countable $2$-colored ultrahomogeneous graphs for which both color classes form disjoint unions of cliques. We identify a new concept, which we call piecewise ultrahomogeneity, as key aspect in such classifications. An ultrahomogeneous graph whose color classes form disjoint unions of cliques is called piecewise ultrahomogeneous if each subgraph induced by a pair of maximal cliques of distinct color is ultrahomogeneous. As explained by Rose (see \cite[Theorem 5.2]{ROS11}), this concept also appears in the dissertation of Jenkinson \cite{JEN06}. We obtain the following characterization of piecewise ultrahomogeneity (see Theorems~\ref{theo:bedingungenuh} and~\ref{theo:hiomittedclassification}):
	
\begin{theoremintro}\label{theo:a}
	Let $G$ be a non-bipartite, countable, 2-colored ultrahomogeneous graph in which the color classes form disjoint unions of cliques and that is not a blow-up. Apart from one degenerate case $F_{2,2}$, the graph $G$ is piecewise ultrahomogeneous if and only if it contains induced subgraphs isomorphic to the graphs $D$ and $\widetilde{D}$ depicted in Figure~\ref{fig:graphs}. 
\end{theoremintro}

We leverage the theorem to completely classify countable $2$-colored ultrahomogeneous graphs in which the color classes form disjoint unions of cliques: 

\begin{theoremintro}\label{theo:b}
	Let $G$ be a countable 2-colored ultrahomogeneous graph in which the color classes form disjoint unions of cliques and that is not a blow-up. Then (after possibly interchanging the colors) exactly one of the following holds:
	
	\begin{enumerate}[(i)]
		\item (Piecewise ultrahomogeneous, Theorem~\ref{theo:twocolors}) Either both color classes in $G$ form an independent set or a single clique, $G$ belongs to a single biparametric family $\{G_{r,b} \colon r,b \in \N \cup \{\aleph_0 \}\}$, or $G$ is isomorphic to the specific graph $F_{2,2}$. 
		\item (Not piecewise ultrahomogeneous, Theorem~\ref{theo:hiomittedclassification}) The graph $G$ belongs to one of two monoparametric families $\{F_{\aleph_0,1}^k \colon k \in \N_{\geq 2}\}$ or $\{F_{\aleph_0,2}^k \colon  k \in \N_{\geq 2}\}$, or it is isomorphic to one of four specific graphs $F_{2,1}$, $F_{\aleph_0,1}$, $F_{\aleph_0,2}$ or~$F_{\aleph_0,\aleph_0}$. 
	\end{enumerate}
\end{theoremintro}

	This paper is organized as follows: Section~\ref{sec:preliminaries} contains preliminary results. In Section~\ref{sec:fraisse}, we recall Fra\"issé's theory and study the structure of minimally omitted subgraphs. In Section~\ref{sec:neighboringrelations}, we introduce the concept of piecewise ultrahomogeneity and prove one implication of Theorem~\ref{theo:a}. In Sections~\ref{sec:omitting} and~\ref{sec:cliquereducible}, we classify graphs that are not piecewise ultrahomogeneous and piecewise ultrahomogeneous, respectively, thereby proving Theorems~\ref{theo:a} and~\ref{theo:b}. We conclude with some final remarks in Section~\ref{sec:conclusion}.
	
\section{Preliminaries}\label{sec:preliminaries}

All graphs in this paper are simple, that is, they neither contain parallel edges nor loops.
Let $G$ be a graph. We denote by $V(G)$ and $E(G)$ the vertex set and the edge set of $G$, respectively.
If $u$ and $v$ are joined by an edge in~$G$, then we write $u \sim _G v$.
We denote the \emph{neighborhood} of a vertex $u$ in $G$ by $N_G(v)$ and simply write $N(v)$ if the ambient graph is clear from the context.
For a subset $S \subseteq V(G)$ we set $N_G^S(v)\coloneqq N_G(v) \cap S$.
If $S \subseteq N_G(v)$, then $v$ \emph{dominates} $S$.
Set $\overline{G}$ to be the edge-complement of~$G$.
The cardinality of a maximum independent set or maximum clique in $G$ is denoted by~$\alpha(G)$ or $\omega(G)$, respectively.
The \emph{lexicographic product} $G \cdot H$ of two graphs $G$ and $H$ is the graph on the vertex set $V(G) \times V(H)$ with $(u_G,u_H)\sim_{G \cdot H}(v_G,v_H)$ if and only if $u_G \sim_G v_G$ or $u_G = v_G$ and $u_H \sim_H v_H$.
The \emph{join} of two disjoint graphs $G$ and $H$ is obtained by adding the edges in $\{uv\colon u \in V(G), v \in V(H)\}$ to the disjoint union of $G$ and $H$.
We denote the complete graph on $n$ vertices by $K_n$, and $P_n$ is a path of order $n$.

\subparagraph{2-colored graphs.}
A \emph{colored graph} is a tuple $(G, \chi_G)$ where $G$ is a graph and $\chi_G$ is a map on $V(G)$.
The sets $\chi_G^{-1}(c)$ for $c \in \chi_G(V(G))$ are the \emph{color classes of $(G, \chi_G)$}.
If $(G, \chi_G)$ has at most $k$ color classes, then $(G, \chi_G)$ is a \emph{$k$-colored graph}.
A 1-colored graph is also called \emph{monochromatic}.
The vertices of induced subgraphs of colored graphs inherit the respective colors.
\medskip

In this article we focus on 2-colored graphs.
We adhere to the convention that the 2-coloring of~$G$ is $\chi_G\colon V(G) \to \{\text{red}, \text{blue}\}$ and set $R_G\coloneqq \chi_G^{-1}(\text{red})$ and $B_G\coloneqq \chi_G^{-1}(\text{blue})$ to be the \emph{red} and \emph{blue} vertices of $G$, respectively.
Two vertices $v,v' \in V(G)$ are \emph{twins} if $\chi_G(v) = \chi_G(v')$ and $N_G(v')\setminus \{v'\} = N_G(v)\setminus \{v\}$.
Edges in $G$ joining vertices of different color are called \emph{cross edges}. We write~$\widetilde{G}$ for the graph obtained from $G$ by complementing the cross edges while maintaining the edges within each color class.
By means of brevity, we often write $G$ instead of $(G, \chi_G)$. We drop the index $G$ whenever the situation is unambiguous.

\subparagraph{(Clique-)Ultrahomogeneity.}
Two colored graphs $G$ and $H$ are \emph{isomorphic} if there exists an \emph{isomorphism} between $G$ and $H$, that is, a bijective color-preserving map $\varphi\colon V(G) \to V(H)$ which satisfies $\varphi(u)\sim \varphi(v)$ if and only if $u \sim v$.
In this case, we write $G \cong H$.
If, additionally, $G = H$, then $\varphi$ is an \emph{automorphism} of $G$.
We write~$\Aut(G)$ for the automorphism group of $G$.
The graph~$G$ is \emph{ultrahomogeneous} if every isomorphism between two finite induced subgraphs of~$G$ extends to an automorphism of $G$.
Note that every color class in a colored ultrahomogeneous graph induces a monochromatic ultrahomogeneous graph.
In order to shorten our notation, we call a graph~$G$ \emph{clique-ultrahomogeneous} (CUH) if $G$ is a countably infinite ultrahomogeneous 2-colored graph where both color classes $R$ and $B$ induce disjoint unions of cliques.
By \cite{LAC80}, the inclusion-wise maximal cliques in $G[R]$ and $G[B]$ are all of the same cardinality $\omega_R\coloneqq \omega(G[R])$ and $\omega_B\coloneqq \omega(G[B])$, respectively.
Setting $\alpha_R \coloneqq \alpha(G[R])$ and $\alpha_B \coloneqq \alpha(G[B])$, we obtain
\[G[R] \cong \overline{K}_{\alpha_R}\cdot K_{\omega_R}\text{ and } G[B]\cong \overline{K}_{\alpha_B}\cdot K_{\omega_B}\text{ with }\max\{\alpha_R,\omega_R,\alpha_B,\omega_B\}=\aleph_0.\]
 We denote the sets of maximal red and blue cliques of~$G$ by $\mR$ and $\mB$, respectively.
Note that $\Aut(G)$ permutes the set $\mR$. Similarly, it permutes $\mB$. From the definition of ultrahomogeneity, we obtain the following statement (also see \cite[Lemma 6.1]{HEI20}):

\begin{lemma}\label{lemma:complement}
	Let $G$ be a 2-colored ultrahomogeneous graph.
	If $H$ is obtained from $G$ by a combination of complementations of the edges within a color class or the cross edges, then $H$ is ultrahomogeneous. 
\end{lemma}

Let $H$ be a 2-colored graph in which one color class is an independent set, say, $R_H$ is an independent set.  We call~$G$ a \emph{blow-up} of $H$ if $G$ is obtained from~$H$ by, for some $i \in \N_{\geq 2} \cup \{\aleph_0\}$, replacing all vertices in this color class by $i$-cliques and joining their vertices to the neighbors of the original vertex in $H$. 
More precisely, $G$ is a blow-up of~$H$ if $G[R] = H[R_H] \cdot K_{\omega_R}$ for some $\omega_R \geq 2$, $G[B] = H[B_H]$, and $(u,v)\sim_G b$ if and only if $u \sim_H b$ for all $u \in V(H)$, $v \in V(K_{\omega_R})$, and $b \in B$.
The following property is easily verified (also see \cite[Lemma 6.2]{HEI20}):

\begin{lemma}\label{lemma:blowupuh}
A blow-up of a graph $H$ is ultrahomogeneous if and only if $H$ is ultrahomogeneous. 
\end{lemma}

We call a CUH graph \emph{basic} if it is not a blow-up and $\min \{\alpha_R, \alpha_B\} \geq 2$.
By complementation inside the color classes and reduction of blow-ups, which preserves ultrahomogeneity (see Lemmas~\ref{lemma:complement} and~\ref{lemma:blowupuh}), we can always pass from any CUH graph to a basic CUH graph. It therefore suffices to consider basic CUH graphs. A 2-colored graph $G$ is \emph{homogeneously connected} if all or none of the possible cross edges in $G$ are present. Concerning the sizes of the color classes, we observe the following:
\begin{lemma}\label{lemma:finitecolorclass}
	Let $G$ be a basic CUH graph. If a color class of $G$ is finite, then $\omega_R = \omega_B = 1$ and $G$ is homogeneously connected. 
\end{lemma}

\begin{proof}
	Without loss of generality, we assume that $|R| \in \mathbb{N}$.
	Since $G$ is countably infinite, we obtain $|B| = \aleph_0$.
	By the pigeonhole principle, there exist distinct blue vertices $b_1$ and $b_2$ with $N^R(b_1) = N^R(b_2)$.
	First assume that $b_1 \sim b_2$. This implies that $\omega_B \geq 2$ and, hence, $G$ is a blow-up, which is a contradiction to~$G$ being basic.
	Hence $b_1 \nsim b_2$ follows.
	If $\omega_B = 1$, then all blue vertices are twins by ultrahomogeneity and, hence, $G$ is homogeneously connected.
	Since $G$ is not a blow-up, we obtain $\omega_R= 1$. 
	If $\omega_B>1$, there exists $b_2' \in B$ with $b_2' \sim b_2$.
	Since $G[\{b_1, b_2\}]\cong G[\{b_1, b_2'\}]$ we obtain that $N^R(b_2) = N^R(b_1) = N^R(b_2')$ and we may apply the above arguments to the adjacent vertices $b_2$ and $b_2'$ to obtain a contradiction.
\end{proof}

However, note that in a basic CUH graph which is not homogeneously connected, either the number or the size of the maximal cliques of a given color can be finite.
\medskip

By~\cite{JEN12}, there exists a unique countably infinite 2-colored ultrahomogeneous graph $G$ with $\omega_R= \omega_B = 1$ which is generic in the following sense: For every $c \in \{\text{red, blue}\}$ and all finite disjoint vertex sets $S, T \subseteq V(G)$ of color $c$, there exists a vertex of color $c' \neq c$ adjacent to all vertices in $S$ and to none of the vertices in $T$. This graph is called the \emph{generic bipartite graph}. 
We frequently make use of the following classification: 

\begin{theorem}[{\cite[Theorem 2.2]{JEN12}}]\label{theo:trussbipartite}
	Let $G$ be a countable 2-colored ultrahomogeneous graph whose color classes form independent sets. Either $G$ is homogeneously connected, the cross edges in $G$ form a perfect matching or its complement, or $G$ is isomorphic to the generic bipartite graph. 
\end{theorem}

Note that the graphs given in Theorem~\ref{theo:trussbipartite} are bipartite.

\section{Fra\"issé limits and omitted subgraphs}\label{sec:fraisse}
In this section, we briefly recall Fra\"issé's theorem and the related terminology. The result as well as further information can be found in standard textbooks on model theory, for example \cite{HOD93}. In the second part of the section, we present a fundamental result on the structure of minimally omitted subgraphs of CUH graphs. 
\medskip

Let $L$ be a countable relational language. An $L$-structure $D$ is \emph{ultrahomogeneous} if every isomorphism between finite substructures extends to an automorphism of~$D$. The \emph{age} $\mathcal{A}_D$ of an $L$-structure~$D$ is the class of all finite structures that are isomorphic to induced substructures of $D$. 
An \emph{amalgamation class} is a class $\mathcal{C}$ of finite $L$-structures which is closed under isomorphism and taking induced substructures, and has the amalgamation property: For $J, A_1, A_2 \in \mathcal{A}$ and embeddings $\iota_i \colon J \to A_i$ ($i = 1,2$), there exists $A \in \mathcal{A}$ and embeddings $\kappa_i \colon A_i \to A$ ($i = 1,2$) such that $\kappa_1 \circ \iota_1 = \kappa_2 \circ \iota_2$ holds. In this situation, $A$ is called an \emph{amalgam} of $A_1$ and~$A_2$.

\begin{theorem}[Fra\"issé]\label{theo:fraisse}
Let $D$ be a countable ultrahomogeneous $L$-structure. Then $\mathcal{A}_D$ is an amalgamation class. Conversely, for every amalgamation class $\mathcal{C}$ of finite $L$-structures, there exists a countable ultrahomogeneous $L$-structure $D$ with $\mathcal{A}_D = \mathcal{C}$, and $D$ is unique up to isomorphism. 
\end{theorem}

In the setting of Theorem~\ref{theo:fraisse}, we call $D$ the \emph{Fra\"issé limit} of $\mathcal{C}$. 
Now we return to the special case of countable 2-colored graphs. By Fra\"issé's theorem, we may shift between countable ultrahomogeneous graphs and amalgamation classes of finite graphs. If $H \in \mathcal{A}_G$, i.e., $H$ is isomorphic to an induced subgraph of $G$, we say that~$H$ is \emph{realized} in $G$. All induced subgraphs of $H$ are then realized in $G$. Conversely, if $H$ is not realized in $G$, we say that $H$ is \emph{omitted} in $G$. In this case, every finite graph $H'$ containing $H$ as induced subgraph is also omitted in $G$. For this reason, it suffices to consider the graphs $H$ which are \emph{minimally omitted} in~$G$: These are the finite graphs $H$ which are omitted in $G$ and for which every proper induced subgraph is realized in $G$. We write $O(G)$ for the set of minimally omitted subgraphs of~$G$.

\begin{Example}
	The countable monochromatic graph $G \cong \overline{K}_s \cdot K_t$ with $s, t \in \mathbb{N} \cup \{\aleph_0 \}$ is ultrahomogeneous.
	The set of minimally omitted subgraphs is 
	\begin{align*}
		O(G) = \begin{cases}
				 \{P_3, K_{t+1}, \overline{K}_{s+1}\} &\text{if $s, t \in \mathbb{N}$}, \\
				 \{P_3, K_{t+1}\} &\text{if $t \in \mathbb{N}, s = \aleph_0$}, \\
				 \{P_3, \overline{K}_{s+1}\} &\text{if $s \in \mathbb{N}, t = \aleph_0$, and} \\
				 \{P_3\} &\text{if $s = \aleph_0, t = \aleph_0$.} 
			   \end{cases}	
	\end{align*}
Omitting $P_3$ forces $G$ to be a disjoint union of cliques. Omitting $K_{t+1}$ or $\overline{K}_{s+1}$ for $s,t \in \N$ restricts the maximal sizes of cliques and independent sets in $G$, respectively.
\end{Example}

As for edge complements, taking ages and omitted sets is compatible with the complementation of the cross edges: 

\begin{lemma}\label{lemma:tildeomitted}
For a countable 2-colored graph $G$, we have $\mathcal{A}_{\widetilde{G}} = \{\widetilde{A} \colon A \in \mathcal{A}_G\}$ and $O(\widetilde{G}) = \{\widetilde{H} \colon H \in O(G)\}$. 
\end{lemma}

The following theorem forms the basis for the arguments in the subsequent sections:

\begin{theorem}\label{theo:twoverticesnoedge}
	Let $G$ be a CUH graph and assume that $H \in O(G)$ is not monochromatic. For every color $c \in \{\text{red, blue}\}$, let $H_c$ be the graph induced by the corresponding color class in $H$. Then one of the following holds: 
	\begin{enumerate}[(i)]
		\item $H_c \cong K_n$ for some $n \geq 3$ and the vertices of $H_c$ are twins in $H$,
		\item $H_c \cong K_2$, or,
		\item $H_c \cong \overline{K}_n$ for some $n \in \mathbb{N}_{\geq 1}$.
	\end{enumerate}
\end{theorem}

	\begin{proof}
		Without loss of generality, let $c$ be blue. The graph $H_{\text{blue}}$ is a disjoint union of cliques since otherwise, a monochromatic $P_3$ is realized in $H_{\text{blue}}$ and, hence, also in $H$. Since monochromatic $P_3$s are omitted in $G$ and $H$ is minimally omitted, it follows that $H$ is a blue $P_3$, which is a contradiction to $H$ not being monochromatic.
		\medskip
		
		We may assume that $|V(H_{\text{blue}})| \geq 3$ since otherwise the claim is trivially satisfied.
		If there exist $b_1,b_2 \in V(H_{\text{blue}})$ with $N^B_H(b_1) \neq N^B_H(b_2)$, then $b_1  \nsim b_2$ since the blue $P_3$ is omitted. By minimality, $H-b_1$ and $H-b_2$ can be embedded into $G$ such that they agree on $H - \{b_1,b_2\}$.
		We may therefore identify $H- \{b_1,b_2\}$ with its image in~$G$ and write $b_1'$ and $b_2'$ for the images of $b_1$ and $b_2$ in the embeddings of $H-b_2$ and $H-b_1$, respectively.
		Since $N^B_H(b_1) \neq N^B_H(b_2)$, the vertices $b_1'$ and $b_2'$ are distinct. Without loss of generality, there exists a blue vertex in	 $v \in N_H(b_1)\setminus N_H(b_2)$.
		We obtain $b_1' \nsim b_2'$ since otherwise, $v b_1'b_2'$ would be an induced blue path.
		This implies that $H-\{b_1,b_2\}$ together with $b_1'$ and $b_2'$ forms an embedding of $H$ into $G$, a contradiction.
		\medskip
		
		It remains to consider the case  that all vertices in $H_{\text{blue}}$ are twins, that is, $H_{\text{blue}}$ is empty or complete. In the first case, the claim is proven, so let $H_{\text{blue}}$ be complete.
		Suppose towards a contradiction that there exist $b_1,b_2 \in V(H_{\text{blue}})$ which are not twins in $H$.
		Embedding $H-b_2$ and $H-b_1$ into $G$ as before, the images $b_1'$ and $b_2'$ of $b_1$ and $b_2$ cannot be equal. However, $b_1$ and $b_2$ have a common blue neighbor $v$, so $b_1'\sim b_2'$ (as $b_1'vb_2´'$ is a blue path otherwise). This means that $H$ is realized in $G$, which is a contradiction. Hence all blue vertices in $H$ are twins.
	\end{proof}
	
\section{Piecewise ultrahomogeneity}\label{sec:neighboringrelations}

We call a CUH graph $G$ \emph{piecewise ultrahomogeneous} if for every $\rC \in \mR$ and $\bC \in \mB$, the graph $G[\rC \cup \bC]$ is ultrahomogeneous. In this section, we show that every basic CUH graph which contains two specific graphs as induced subgraphs is piecewise ultrahomogeneous. 
\medskip

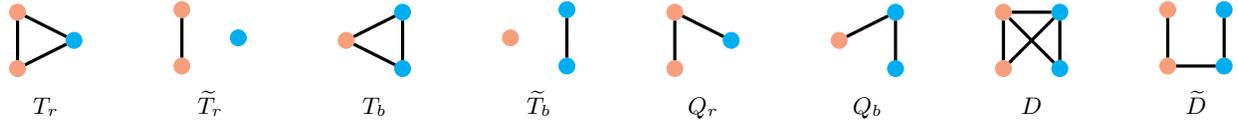
\begin{figure}
	\centering
	
	\begin{subfigure}{1.2cm}
		\centering
		\begin{tikzpicture}[scale=.5, yscale=1.5]
			\node[\myred] (r1) at (0, 1) {};
			\node[\myred] (r2) at (0, 0) {};
			\node[\myblue] (b1) at (1.5, .5) {};
			
			\draw[edge] (r1) to (b1) to (r2) to (r1);
		\end{tikzpicture}
		\subcaption*{$T_r$}
	\end{subfigure}
	\hspace{.75cm}
	\begin{subfigure}{1.2cm}
		\centering
		\begin{tikzpicture}[scale=.5, yscale=1.5]
			\node[\myred] (r1) at (0, 1) {};
			\node[\myred] (r2) at (0, 0) {};
			\node[\myblue] (b1) at (1.5, .5) {};
			
			\draw[edge] (r1) edge (r2);
		\end{tikzpicture}
		\subcaption*{$\widetilde{T}_r$}
	\end{subfigure}
	\hspace{.75cm}
	\begin{subfigure}{1.2cm}
		\centering
		\begin{tikzpicture}[scale=.5, yscale=1.5]
			\node[\myred] (r1) at (0, .5) {};
			\node[\myblue] (b1) at (1.5, 1) {};
			\node[\myblue] (b2) at (1.5, 0) {};
			
			\draw[edge] (r1) to (b1) to (b2) to (r1);
		\end{tikzpicture}
		\subcaption*{$T_b$}
	\end{subfigure}
	\hspace{.75cm}
	\begin{subfigure}{1.2cm}
		\centering
		\begin{tikzpicture}[scale=.5, yscale=1.5]
			\node[\myred] (r1) at (0, .5) {};
			\node[\myblue] (b1) at (1.5, 1) {};
			\node[\myblue] (b2) at (1.5, 0) {};
			
			\draw[edge] (b1) to (b2);
		\end{tikzpicture}
		\subcaption*{$\widetilde{T}_b$}
	\end{subfigure}
	\hspace{.75cm}
	\begin{subfigure}{1.2cm}
		\centering
		\begin{tikzpicture}[scale=.5, yscale=1.5]
			\node[\myred] (r1) at (0, 1) {};
			\node[\myred] (r2) at (0, 0) {};
			\node[\myblue] (b1) at (1.5, 0.5) {};
			
			\draw[edge] (b1)--(r1)--(r2);
		\end{tikzpicture}
		\subcaption*{$Q_r$}
	\end{subfigure}
	\hspace{.75cm}
	\begin{subfigure}{1.2cm}
		\centering
		\begin{tikzpicture}[scale=.5, yscale=1.5]
			\node[\myred] (r1) at (0, 0.5) {};
			\node[\myblue] (b1) at (1.5, 1) {};
			\node[\myblue] (b2) at (1.5, 0) {};
			
			\draw[edge] (r1)--(b1)--(b2);
		\end{tikzpicture}
		\subcaption*{$Q_b$}
	\end{subfigure}
	\hspace{.75cm}
	\begin{subfigure}{1.2cm}
		\centering
		\begin{tikzpicture}[scale=.5, yscale=1.5]
			\node[\myred] (r1) at (0, 1) {};
			\node[\myred] (r2) at (0, 0) {};
			\node[\myblue] (b1) at (1.5, 1) {};
			\node[\myblue] (b2) at (1.5, 0) {};
			
			\draw[edge]
			(r1) edge (r2)
			edge (b1)
			edge (b2)
			(r2) edge (b1)
			(b1) edge (b2);
		\end{tikzpicture}
		\subcaption*{$D$}
	\end{subfigure}
	\hspace{.75cm}
	\begin{subfigure}{1.2cm}
		\centering
		\begin{tikzpicture}[scale=.5, yscale=1.5]
			\node[\myred] (r1) at (0, 1) {};
			\node[\myred] (r2) at (0, 0) {};
			\node[\myblue] (b1) at (1.5, 1) {};
			\node[\myblue] (b2) at (1.5, 0) {};
			
			\draw[edge] (r1) to (r2) to (b2) to (b1);
		\end{tikzpicture}
		\subcaption*{$\widetilde{D}$}
	\end{subfigure}

	\caption{The graphs in $\mathcal{T} \cup \{Q_r, Q_b, D,\widetilde{D}\}$}\label{fig:graphs}	
\end{figure}

Let $T_r$ and $T_b$ be the triangles containing a single blue vertex and a single red vertex, respectively. We set $\mathcal{T} = \{T_r, \widetilde{T}_r, T_b, \widetilde{T}_b\}$. Let $Q_r$ and $Q_b$ be the graphs obtained from $T_r$ and $T_b$, respectively, by deleting one cross edge. Moreover, let $D$ be the graph arising from a complete graph on two red and two blue vertices by deleting one cross edge (that is, a 2-colored diamond in which no two vertices of the same color are of the same degree). The graphs in $\mathcal{T}$ as well as $Q_r$, $Q_b$, $D$ and $\widetilde{D}$ are depicted in Figure~\ref{fig:graphs}.

\medskip

The aim of this section is the proof of the following result:

\begin{theorem}\label{theo:bedingungenuh}
	Let $G$ be a basic CUH graph. If $D$ and $\widetilde{D}$ are realized in $G$, then $G$ is piecewise ultrahomogeneous. 
\end{theorem}

To increase the readability, we introduce the following convention: Until the end of this section, we assume that $G$ is a basic CUH graph in which $D$ and $\widetilde{D}$ are realized. Moreover, we formulate the statements of Remark~\ref{rem:realizedgraphs} as well as Lemmas~\ref{lemma:nonadjointneighbors}, \ref{lemma:onebluevertex} and \ref{lemma:kneighbors} for the red vertices, but the analogous results hold for the blue color class. 

\begin{Remark}\label{rem:realizedgraphs}
$\null$
\begin{enumerate}[(i)]
\item We have $\omega_R, \omega_B \geq 2$.
\item Fix a maximal blue clique $\bC \in \mB$, and let $r_1, r_2 \in R$ with $r_1 \sim r_2$.
	We claim that in $G[\{r_1, r_2\} \cup \bC]$ all of the following graphs are realized: $T_r$, $\widetilde{T}_r$, and $Q_r$. Moreover, $Q_r$ can be realized in both ways -- either of the vertices $r_1, r_2$ can correspond to the degree-2-vertex of $Q_r$.
    To see this, fix $b \in \bC$.
    First assume that $b \in N(r_1)\setminus N(r_2)$.
	Using ultrahomogeneity and the assumption that $D$ is realized in $G$, we find that $T_r$ is realized in $G[\{r_1, r_2\} \cup \bC]$. Arguing similarly for $\widetilde{D}$, we obtain that $\widetilde{T}_r$ is realized in $G[\{r_1, r_2\} \cup \bC]$.
	We combine this with $\widetilde{D}$ being realized in $G$ to obtain that $(N(r_2)\setminus N(r_1))\cap \bC \neq \emptyset$, that is, $Q_r$ is realized in both ways as desired.
	 Now if $b\in N(r_1)\cap N(r_2)$, we find a vertex $b' \in \bC$ which is adjacent to precisely one vertex in $\{r_1, r_2\}$ by using that $D$ is realized in $G$. Similarly, we argue if $b$ is a joint non-neighbor of $r_1$ and $r_2$. 
\end{enumerate}
\end{Remark}

	\begin{lemma}\label{lemma:nonadjointneighbors}
Let $C \subseteq R$ be a finite clique.
If the vertices in $C$ have a joint blue neighbor, then they also have a pair of non-adjacent joint blue neighbors.
	\end{lemma}
	
	\begin{proof}
We prove the claim by induction on $|C|$. For $|C| \in \{1,2\}$, the claim follows from Remark~\ref{rem:realizedgraphs}. Let $|C| \geq 3$.
By assumption, there is a vertex $b^{\star}$ in some blue clique $\bC \in \mB$ with $C \subseteq N_G(b^{\star})$.
By Remark~\ref{rem:realizedgraphs} there is a vertex $v \in \bC$ with $1 \leq |N_G(v)\cap C| \leq |C|-1$.
Consider a 2-colored graph $H$ whose red vertices form a $|C|$-clique and whose blue vertices form an independent 2-set $\{b_1, b_2\}$ with $N_H(b_1) = R_H$ and $|N_H(b_2)| = |N_G(v)\cap C|$.
If $H$ is realized in $G$, then we may see $H$ as an induced subgraph of $G$.
Mapping the red vertices of $H$ to $C$  and $b_2$ to $v$ yields a partial isomorphism of $G$, which extends to an automorphism of $G$ mapping $b_1$ to some vertex $b_1' \notin \bC$.
We obtain that $b^{\star}$ and $b_1'$ satisfy the claim.
\medskip

If $H$ is omitted in $G$, then there exists an induced subgraph $H'$ of $H$ with $H' \in O(G)$.
By construction, $H-b_1$ and $H-b_2$ are realized in $G$ and, hence, $\{b_1, b_2\} \subseteq  V(H')$.
We obtain from Remark~\ref{rem:realizedgraphs} that $|R_{H'}| \geq 3$ and, hence, all vertices in $R_{H'}$ are twins (see Theorem~\ref{theo:twoverticesnoedge}).
In particular $|R_{H'}| \leq |C|-1$ since the vertices of $C$ are not twins in $H$.
For $i \in \{1,2\}$ either $R_{H'} = N_{H'}(b_i)$ or $R_{H'} \cap N_{H'}(b_i) = \emptyset$.
By construction we have $R_{H'} = N_{H'}(b_1)$.
If $R_{H'} = N_{H'}(b_2)$, then $H' \in O(G)$ yields a contradiction to the induction hypothesis.
Hence $R_{H'} \cap N_{H'}(b_2) = \emptyset$.
Fix a red $|R_{H'}|$-clique $\hat{C}$ in $G$.
Since $H'$ is minimally omitted, there is a blue vertex $\hat{b}$ with $\hat{C}\cap N_G(\hat{b}) = \emptyset$.
By induction, there exist non-adjacent blue vertices $\hat{b}_1$ and $\hat{b}_2$ such that $\hat{C} \subseteq N_G(\hat{b}_i)$ for $i \in \{1,2\}$.
Since $H'$ is omitted, we have $\hat{b}_i \sim \hat{b}$ for $i \in \{1,2\}$. But then $\hat{b}_1\hat{b}\hat{b}_2$ is a monochromatic $P_3$ realized in~$G$, a contradiction.
	\end{proof}

\begin{lemma}\label{lemma:onebluevertex}
Let $H$ be the join of a blue $K_1$ and a finite red $K_ {\ell}$ with $\ell \leq \omega_R$.
Both graphs $H$ and $\widetilde{H}$ are realized in $G$. 
\end{lemma}
	
	\begin{proof}
	 We proceed by induction on~$\ell$.
 	For $\ell \in \{1,2\}$ the claim follows with Remark~\ref{rem:realizedgraphs}.
 	Let $\ell \geq 3$. Suppose that~$H$ is omitted in $G$ and let $H' \in O(G)$ be an induced subgraph of $H$.
 	Observe that $H'$ contains a blue vertex since a red $K_{\ell}$ is realized in $G$.
  Set $\ell' \coloneqq |R_{H'}|$.
  By ultrahomogeneity and since $H' \in O(G)$ it follows that every blue vertex has precisely $\ell'-1$ neighbors in every maximal red clique.
  Fix $\rC \in \mR$ and $b_1,b_2 \in B$ with $b_1\sim b_2$.
  Remark~\ref{rem:realizedgraphs} yields $N_G^{\rC}(b_1) \neq N_G^{\rC}(b_2)$. Moreover, there exists a vertex $r \in \rC \cap N(b_1)\cap N(b_2)$. 
  By Lemma~\ref{lemma:nonadjointneighbors}, there exists vertex $b_N$ dominating $N_G^{\rC}(b_1)$ in a maximal blue clique $\bC' \neq \bC$. 
  Since $|N^{\rC}_G(b_N)| = \ell'-1$, we obtain $N_G^{\rC}(b_1) = N_G^{\rC}(b_N)$.
  Fixing $b_N$ and $r$ and exchanging $b_1$ and $b_2$ defines a partial isomorphism~$\varphi$ of $G$.
  Let $\hat{\varphi} \in \Aut(G)$ be an extension of~$\varphi$.
  Due to $\varphi(r) = r$, the map $\hat{\varphi}$ fixes $\rC$ setwise and, hence, also $N_G^{\rC}(b_1) = N_G^{\rC}(b_N)$ is fixed setwise.
  On the other hand, $\hat{\varphi}(b_1) =b_2$ implies $\hat{\varphi}(N_G^{\rC}(b_1)) = N_G^{\rC}(b_2)$, which contradicts $N_G^{\rC}(b_1) \neq N_G^{\rC}(b_2)$.
  The case of~$\widetilde{H}$ can be treated similarly.
	\end{proof}

\begin{lemma}\label{lemma:kneighbors}
	Every graph consisting of a blue $K_1$, a finite red $K_ {\ell}$ with $\ell \leq \omega_R$, and arbitrary cross edges is realized in~$G$.
\end{lemma}

\begin{proof}
Let $H$ be such a graph. 
We denote the number of cross edges in $H$ by $k \in \{0, \ldots, \ell\}$.
By Lemma~\ref{lemma:onebluevertex} we may assume $k \in \{1, \ldots, \ell-1\}$. By Lemma~\ref{lemma:onebluevertex}, $G$ contains a complete induced subgraph $K$ consisting of a red $k$-clique~$C$ and a blue vertex $b$. Let $\rC \in \mR$ with $C \subseteq \rC$. By assumption, $b$ has a non-neighbor in $\rC$. By Lemma~\ref{lemma:onebluevertex}, a disjoint union of a red $K_{\ell-k}$ and a blue $K_1$ is realized as induced subgraph of $G$. By ultrahomogeneity, there exist distinct vertices $w_1, \ldots, w_{\ell-k} \in \rC$ to which $b$ is non-adjacent. But then $H \cong G[V(K) \cup \{w_1, \ldots, w_{\ell-k}\}]$ and, hence, $H$ is realized in $G$.
\end{proof}

\begin{lemma}\label{lemma:neighborseveryclique}
Fix $\bC \in \mB$ and a finite red $\ell$-clique $C$ for $\ell \leq \omega_R$. For every $S \subseteq C$, there exists $v \in \bC$ with $N_G^C(v) = S$. 
\end{lemma}
	
\begin{proof}
We first show that for every $k \in \{0, \ldots, \ell\}$, there exists $v_k \in \bC$ with $|N_G^C(v_k)| = k$. To this end, we proceed by induction on $\ell$. For $\ell \in \{1,2\}$ the claim follows by Remark~\ref{rem:realizedgraphs}. Now let $\ell \geq 3$. By Remark~\ref{rem:realizedgraphs}, there exists a vertex $w \in \bC$ with $|N_G^C(w)| = j$ for some $j \in \{1, \ldots, \ell-1\}$. Now assume that for some $k \in \{0, \ldots, \ell\}$, none of the vertices in $\bC$ has precisely $k$ neighbors in $C$. Consider the finite graph $H$ with $V(H) = \{h_x \colon x \in C \cup \{w, v\}\}$ such that $G[C \cup \{w\}] \to (H-h_{v}), \ x \mapsto h_x$ is an isomorphism. The blue vertex $h_{v}$ is adjacent to $h_{w}$ and has precisely $k$ neighbors in $\{h_x \colon x \in C\}$. As in the proof of Lemma~\ref{lemma:nonadjointneighbors}, it follows that $H$ is omitted in $G$. Let $H'$ be an induced subgraph of $H$ that is minimally omitted in $G$. By the induction hypothesis, $H'$ contains $h_x$ for all $x \in C$. The graphs $H-h_{w}$ and $H-h_{v}$ are realized in $G$ by Lemma~\ref{lemma:kneighbors}. Hence~$H$ is minimally omitted. Due to $j \in \{1, \ldots, \ell-1\}$, the red vertices in $H'$ are not twins. This is a contradiction to Theorem~\ref{theo:twoverticesnoedge}.
\medskip

Now let $S \subseteq C$ be an arbitrary subset. By the first part of this proof, there exists a vertex $v' \in \bC$ such that $S' \coloneqq N_G^C(v')$ has size $|S|$. Moreover, there exists vertex $b \in \bC$ dominating $C$. Now consider the partial isomorphism $\varphi$ of $G$ obtained by bijectively mapping $S'$ to $S$ and $C \setminus S'$ to $C \setminus S$ while fixing $b$. Let $\hat{\varphi} \in \Aut(G)$ be an extension of~$\varphi$ to~$G$. Then $v \coloneqq \hat{\varphi}(v') \in \bC$ is a vertex with $N_G^C(v) = S$.
\end{proof}
With these results, we prove Theorem~\ref{theo:bedingungenuh}:

\begin{proof}[Proof of Theorem~\ref{theo:bedingungenuh}]
Let $\rC \in \mR$ and $\bC \in \mB$. We show that $G[\rC \cup \bC]$ is the complement of the generic bipartite graph. Consider finite disjoint subsets $S, T \subseteq \rC$.
Since $S \cup T$ is a red clique of size $|S| + |T|$, there exists a vertex $v \in \bC$ with $N_G^{S \cup T}(v) = S$ (see Lemma~\ref{lemma:neighborseveryclique}). In other words, $v$ is adjacent to all vertices in $S$ and to none of the vertices in~$T$.
For the blue color class, one can argue similarly. Hence the complement of $G[\rC \cup\bC]$ is the generic bipartite graph.
\end{proof}

\section{\texorpdfstring{Graphs omitting $D$ or $\widetilde{D}$}{Graphs omitting D or Dtilde}}\label{sec:omitting}
In this section, we classify the basic CUH graphs which omit $D$ or $\widetilde{D}$. We prove that such a graph is either isomorphic to one of the graphs given in Theorem~\ref{theo:trussbipartite}, it belongs to one of two monoparametric families, or it is isomorphic to one of five specific graphs (see Theorem~\ref{theo:hiomittedclassification}). In Section~\ref{sec:minimallyomitted}, we determine the structure of possible minimally omitted subgraphs. In Section~\ref{sec:classification}, we prove our classification result. 

\subsection{Structure of minimally omitted subgraphs}\label{sec:minimallyomitted}

Recall the definition of the graphs in $\mathcal{T}$ as well as $D$ and $\widetilde{D}$ from Section~\ref{sec:neighboringrelations}. The aim of this subsection is the proof of the following theorem: 

\begin{theorem}\label{theo:omittedgraphs}
Let $G$ be a basic CUH graph in which $D$ or $\widetilde{D}$ is omitted, and which is not isomorphic to one of the graphs in Theorem~\ref{theo:trussbipartite}. Then every non-monochromatic graph in $O(G)$ is contained in $\mathcal{T} \cup \{D, \widetilde{D}\}$. Moreover, we have $H \in O(G)$ if and only if $\widetilde{H} \in O(G)$ for $H \in \mathcal{T} \cup \{D, \widetilde{D}\}$. In particular, $G = \widetilde{G}$.
\end{theorem}

Throughout this section, we assume that $G$ is a basic CUH graph that omits $D$ and in which one of the color classes does not form an independent set. The case that $G$ omits $\widetilde{D}$ will be considered at the end of this section. We proceed as follows: In Section~\ref{sec:tilde}, we investigate the structure of the graphs in $O(G)$, with a particular focus on cross edge complements. In Section~\ref{sec:partitions}, we study certain partitions induced by the neighboring relations in~$G$. These results are used to prove Theorem~\ref{theo:omittedgraphs} (see Section~\ref{sec:proofomitted}). For the sake of readability, we formulate several results only for the red color class, but the analogous version for the blue vertices holds as well. 

\subsubsection{Minimally omitted subgraphs and cross edge complements}\label{sec:tilde}

Throughout this subsection, we assume that $G$ is a basic CUH graph that omits $D$
and in which one of the color classes does not form an independent set. We now study the structure of the graphs in $O(G)$. 

\begin{Remark}\label{rem:o3}
If $Q_r \in O(G)$, then all pairs of adjacent red vertices have the same blue neighbors and, hence,~$G$ is a blow-up.
	This is a contradiction to $G$ being basic. We thus obtain $Q_r \notin O(G)$, and, similarly, $Q_b \notin O(G)$. 
\end{Remark}

Hence, since $G$ is basic and there exists an induced subgraph of $D$ which is minimally omitted in~$G$, one of the following cases arises:
\begin{enumerate}[(O1)]
	\item Precisely one of the color classes in $G$ forms an independent set. 
	\item The graph $G$ minimally omits one of the graphs $T_r$ or $T_b$.
	\item The graph $G$ minimally omits $D$. 
\end{enumerate}

%Using ultrahomogeneity together with the fact that $G$ is not a blow-up, it is easy to see that the following holds:

\begin{lemma}\label{lemma:localstructure}
If $\omega_R> 1 $ holds, then every blue vertex has both neighbors and non-neighbors in every maximal red clique. In particular, we obtain
$G[\rC \cup \bC] \cong G[\rC' \cup \bC']$ for all $\rC,\rC' \in \mR$ and $\bC,\bC' \in \mB$.
\end{lemma}
\begin{proof}
	Since $G$ is not a blow-up and $\omega_R>1$ holds, we obtain $N^{\rC}(b) \neq \emptyset \neq \rC \setminus N^{\rC}(b)$ for every $b \in B$ and every $\rC \in \mR$.
	In particular, for all $\rC,\rC' \in \mR$ and $\bC,\bC' \in \mB$ there exists a partial isomorphism $\varphi$ mapping an edge in $G[\rC \cup \bC]$ to an edge in $G[\rC' \cup \bC']$.
	By the ultrahomogeneity of $G$, the isomorphism $\varphi$ extends to an automorphism $\hat{\varphi} \in \Aut(G)$.
	Restricting $\hat{\varphi}$ to $\rC \cup \bC$ yields the desired isomorphism.
\end{proof}

We now consider the case (O2).

\begin{lemma}\label{lemma:nojointneighbors}
We have $T_r \in O(G)$ if and only if $\widetilde{T}_r \in O(G)$ holds, and this is the case precisely if $\omega_R = 2$ holds. 
\end{lemma}

\begin{proof}
	For $\omega_R = 1$, the claim trivially holds, so assume $\omega_R \geq 2$.
	Suppose that $T_r \in O(G)$ holds, that is, every blue vertex has at most one neighbor in each red clique.
	Let~$r \in R$ and assume that $r$ has non-adjacent blue neighbors $b_1$ and $b_2$.
	Denote by~$\rC$ the maximal red clique containing~$r$ and consider a red vertex $r' \in \rC\setminus \{r\}$. By ultrahomogeneity, $r'$ has non-adjacent blue neighbors, so there exists $b \in N_G^B(r')$ with $b \nsim b_1$. 
    If there exists a vertex $v \in \rC \setminus (N(b_2) \cup N(b))$, then fixing $b_1$ and $v$, and exchanging $b_2$ and $b$ defines a partial isomorphism $\varphi$ of $G$. Consider an extension $\hat{\varphi} \in \Aut(G)$ of $\varphi$ to~$G$. Since $v$ is fixed, $\hat{\varphi}$ fixes $N^B(v)$ setwise. Since $b_1$ is fixed, we have $\hat{\varphi}(r) = r$. But on the other hand, we have $\hat{\varphi}(r) = r'$ as $b_2$ is mapped to $b$. This is a contradiction. Hence $\rC \setminus (N(b_2) \cup N(b)) = \emptyset$.
    This yields $\omega_R = |\rC| = 2$, and both $b_2$ and~$b$ are adjacent to precisely one of the two vertices in $\rC$. 
	It remains to consider the case that $N^B(r)$ is contained in a single maximal blue clique $\bC \in \mB$. By Lemma~\ref{lemma:localstructure}, this forces $\omega_B = 1$,
	%Using ultrahomogeneity together with $\alpha_B\geq 2$, it is easily verified that $N^B(r) = \bC$. If $\omega_B > 1$ holds, then $G$ is a blow-up, a contradiction. Otherwise, 
	so $|N_G^B(r)| = 1$. It is easy to see that this is impossible. Altogether, we obtain that $T_r \in O(G)$ implies $\omega_R=2$. Replacing $b_1$ and $b_2$ in the above argument by non-adjacent blue non-neighbors of~$r$ and proceeding analogously shows that $\widetilde{T}_r \in O(G)$ implies $\omega_R = 2$.
	\medskip
	
	Conversely, assume $\omega_R = 2$.
	If there is a blue vertex adjacent to both or none of the vertices of a maximal red clique, then ultrahomogeneity implies that $G$ is a blow-up, a contradiction. Hence $T_r$ and $\widetilde{T}_r$ are omitted in $G$.
	Observe that neither a red $K_2$ is omitted in $G$ (due to $\omega_R = 2$) nor a 2-colored $K_2$ or $\overline{K}_2$ is omitted in $G$ (since $G$ is not a blow-up).
	Thus, $\{T_r, \widetilde{T}_r\} \subseteq O(G)$.
\end{proof}

Now we study the case (O3). There, $T_r$ and $T_b$ are realized in $G$. By Lemma~\ref{lemma:nojointneighbors}, all graphs in $\mathcal{T}$ are realized in $G$.

\begin{lemma}\label{lemma:h1h2}
	If $D \in O(G)$, then $\widetilde{D} \in O(G)$.
\end{lemma}

\begin{proof}
Suppose towards a contradiction that $\widetilde{D}$ is realized in $G$.
Let $r_1, r_2 \in R$ with $r_1 \sim r_2$. Since $T_r$ is realized, there exists a vertex $b \in N^B(r_1) \cap N^B(r_2)$.
Let $\bC \in \mB$ be the maximal blue clique containing $b$. Since $G$ is not a blow-up, there exists a vertex $b' \in \bC \setminus (N^B(r_1) \cap N^B(r_2))$. Then $b' \in \bC \setminus (N^B(r_1) \cup N^B(r_2))$ since  $D \in O(G)$. Hence $N^{\bC}(r_1) = N^{\bC}(r_2)$ follows. By ultrahomogeneity, all vertices in the maximal red clique containing $r_1$ and~$r_2$ have the same neighbors in $M_B$. This is a contradiction to Lemma~\ref{lemma:localstructure}. Thus $\widetilde{D}$ is omitted, and, hence, minimally omitted, in $G$. 
\end{proof}

\subsubsection{Partitions}\label{sec:partitions}
As before, we assume that $G$ is a basic CUH graph omitting $D$ in which one of the color classes does not form an independent set. We study certain partitions of the maximal monochromatic cliques which are induced by neighboring relations. Using these, we deduce that $\omega_R, \omega_B \in \{1,2,\aleph_0\}$ holds.

\begin{lemma}\label{lem:partition}
	Let $\rC \in \mR$ and $\bC \in \mB$.
	If $\omega_R >1$, there is a partition $\rC = R_1 \dot{\cup} R_2$ into non-empty parts such that $N^{\rC}(b) \in \{R_1, R_2\}$ holds for every $b \in \bC$. 
	If, additionally, $\omega_B >1$ holds, there exist partitions $\rC = R_1 \dot{\cup} R_2$ and $\bC = B_1 \dot{\cup} B_2$ into non-empty parts such that $(R_1 \times B_1) \cup (R_2 \times B_2)$ is the set of cross edges in~$G[\rC \cup \bC]$.
\end{lemma}

\begin{proof}
Let $\omega_R> 1$. If $T_r$ or $\widetilde{T}_r$ are omitted in $G$, the claim follows by Lemma~\ref{lemma:nojointneighbors}.
Moreover, the statement follows immediately for $\omega_B = 1$ since $G$ is not a blow-up.
From now on, we assume $\omega_B >1$ and that $T_r$ and $\widetilde{T}_r$ are realized in $G$.
\medskip

If $T_b$ is omitted in~$G$, then $T_b \in O(G)$ since $\omega_B >1$ and $G$ is not a blow-up.
With Lemma~\ref{lemma:nojointneighbors} we obtain $\omega_B = 2$. Let $b_1, b_2 \in \bC$ be distinct vertices. Every vertex in $\rC$ is adjacent to precisely one vertex in $\{b_1, b_2\}$. Defining $R_1 \coloneqq N^{\rC}(b_1)$ and $R_2 \coloneqq N^{\rC}(b_2)$ yields the desired partition.
We proceed analogously if $\widetilde{T}_b$ is omitted in~$G$.
\medskip

Now assume that the graphs in $\mathcal{T}$ are realized in $G$. By Lemma~\ref{lemma:h1h2}, we have $D, \widetilde{D} \in O(G)$.
Consider the equivalence relation on $\rC$ which contains $(r_1, r_2)$ if and only if $N^{\bC}(r_1) = N^{\bC}(r_2)$.
Since $D$ and $\widetilde{D}$ are minimally omitted in $G$, not all vertices in $\rC$ have the same neighbors in $\bC$, so there are at least two equivalence classes.
On the other hand, let $r_1$ and $r_2$ be in different equivalence classes. Then every vertex in $\bC$ is joined to $r_1$ precisely if it is not joined to $r_2$ since  $D, \widetilde{D} \in O(G)$. Hence there are exactly two equivalence classes. This means that $\rC$ is partitioned into two subsets and every vertex in $\bC$ is adjacent to precisely one of these. The second claim is a simple consequence of the first part. 
\end{proof}

Assume $\omega_R >1$. Let $\rC \in \mR$ and $\bC,\bC' \in \mB$. We say that $\bC$ and $\bC'$ \emph{induce the same partition in $\rC$} if there exist $b \in \bC$ and $b' \in \bC'$ such that $N^{\rC}(b) \in \{N^{\rC}(b'), \rC \setminus N^{\rC}(b') \}$ (if $\omega_B >1$ holds, then we can require $N^{\rC}(b) = N^{\rC}(b')$).

\begin{Remark}\label{rem:samepartition}
Assume $\omega_R >1$ and let $\bC,\bC' \in \mB$ be distinct maximal cliques inducing the same partition in $\rC \in \mR$.
By ultrahomogeneity all maximal blue cliques induce the same partition in $\rC$. If $\omega_R > 2$, then there exist pairs of vertices in $\rC$ that have a joint blue neighbor, and others that do not. This is a contradiction and hence $\omega_R = 2$ follows. 
\end{Remark}

\begin{theorem}\label{theo:sizeclique}\label{theo:finitecliques}
We have $\omega_R \in \{1,2,\aleph_0\}$.
If $\omega_R = \aleph_0$, then $\alpha_B= \aleph_0$ and we obtain $|N_G^{\rC}(b)| = |\rC \setminus N_G^{\rC}(b)| = \aleph_0$ for every $b \in B$ and every $\rC \in \mR$.
\end{theorem}

\begin{proof}
Let $\rC \in \mR$.
Suppose towards a contradiction that $\omega_R \in \mathbb{N}_{\geq 3}$.
With Remark~\ref{rem:samepartition} we obtain $\alpha_B \in \mathbb{N}$.
Let $b \in B$ and set $k \coloneqq |N^{\rC}(b)|$.
For every partition of $\rC$ into sets of size $k$ and $\omega_R-k$, we find precisely one maximal blue clique inducing this partition. By ultrahomogeneity, this holds for all maximal red cliques. Hence every permutation of a maximal red clique induces a unique permutation of the maximal blue cliques.
Pointwise fixing $\rC' \in \mR\setminus \{\rC\}$ and cyclically permuting $\rC$ therefore defines a partial isomorphism of $G$ that cannot be extended to an automorphism of $G$. 
\medskip

Now suppose $\omega_R = \aleph_0$. We show $|N^{\rC}(b)| = \aleph_0$ for all $b \in B$. If $\omega_B > 1$, this follows by ultrahomogeneity and Lemma~\ref{lem:partition}. Now assume $\omega_B = 1$, and suppose $|N^{\rC}(b)| = k$ for some $k \in \N$. By Lemma~\ref{lemma:nojointneighbors}, we have $k \geq 2$. By ultrahomogeneity, every blue vertex has precisely $k$ neighbors in every maximal red clique. Moreover, for every subset $S \subseteq \rC$ of size $k$, there exists $b_S \in B$ with $N_G^{\rC}(b_S) = S$. Fix such a set $S$ and let $S_1, S_2  \subseteq \rC$ be sets of size $k$ with $|S \cap S_1| = k-1$ and $|S \cap S_2| = k-2$.
Since $\rC$ is infinite, there exists a vertex $r \in \rC \setminus (N^{\rC}(b) \cup N^{\rC}(b_{S_1}) \cup N^{\rC}(b_{S_2}) )$.
Fixing~$r$ and mapping $b_S$ to $b_{S_1}$, $b_{S_1}$ to $b_{S_2}$, and $b_{S_2}$ to $b_S$ defines a partial isomorphism of~$G$. Its extension to $G$ maps $S_2 \cap S$ to $S \cap S_1$, which is a contradiction since these sets do not have the same cardinality. Similarly, one shows $|M_R \setminus N^{\rC}(b)| = \aleph_0$. 
\medskip

It remains to show $\alpha_B= \aleph_0$. This directly follows if $\omega_B \in \{1,2\}$. Assume $\omega_B = \aleph_0$ and fix $\rC \in \mR$. If $\alpha_B$ is finite, there exist distinct vertices in $\rC$ which have the same blue neighbors. On the other hand, not all vertices in $\rC$ are twins since $G$ is not a blow-up. This contradicts the ultrahomogeneity of $G$.
\end{proof}

\subsubsection{\texorpdfstring{Proof of Theorem~\ref{theo:omittedgraphs}}{Proof of Theorem 5.1}}\label{sec:proofomitted}
In this section, we prove Theorem~\ref{theo:omittedgraphs}. We mainly show that every non-monochromatic graph in $O(G)$ is contained in $\mathcal{T} \cup \{D, \widetilde{D}\}$. To this end, we use the results from Section~\ref{sec:tilde} together with the fact that $\omega_R, \omega_B \in \{1,2,\aleph_0\}$ by Theorem~\ref{theo:sizeclique}. We begin with a technical lemma:

\begin{lemma}\label{lemma:clique+indepset}
Let $H$ be the join of a red $K_r$ and a blue $\overline{K}_k$ for $r,k \in \N$ with $k \leq \alpha_B$. If  $\omega_R = \aleph_0$, then $H$ and $\widetilde{H}$ are realized in $G$.  
\end{lemma}

\begin{proof}
	Assume $\omega_B > 1$. Fix $\rC \in \mR$ and distinct maximal blue cliques $\bC^1, \ldots, \bC^k \in \mB$. Consider the partition of $\rC$ obtained by intersecting the partitions induced by $\bC^1, \ldots, \bC^k$. Since $\rC$ is infinite, at least one part $P$ is infinite. The vertices in $P$ have the same neighbors in $\bC^1 \cup \ldots \cup \bC^k$. Choose $b_i \in N_G^{\bC^i}(P)$ for $i \in \{1, \ldots, k\}$, and let $P' \subseteq P$ be a set of size~$r$. Then $G[\{b_1, \ldots, b_k\} \cup P']$ is isomorphic to $H$.
	\medskip
	
	Now assume $\omega_B= 1$.
	For a contradiction, suppose that $G$ omits $H$.
	We may assume without loss of generality that $H \in O(G)$ by replacing $H$ by a smallest omitted subgraph induced in $H$ (which is again a join or a disjoint union of a complete red graph and an edgeless blue graph). 
	For a contradiction, suppose $H \in O(G)$. Fix an $r$-clique~$C \subseteq R$ and let $b_1, \ldots, b_{k-1}$ be the blue vertices dominating $C$. Let $\rC^C \in \mR$ be the maximal red clique containing $C$ and consider $v \in \rC^C \setminus N(b_1)$. There are at most $k-2$ blue vertices dominating the red clique $C \cup \{v\}$. By ultrahomogeneity, this holds for every red $(r+1)$-clique. Continuing this way, we find a red $r'$-clique for some $r' \geq r$ which is not dominated by a blue vertex. This is a contradiction to Theorem~\ref{theo:sizeclique}. The proof for $\widetilde{H}$ is similar.
\end{proof}

\begin{lemma}\label{lemma:independentseth}
Assume that $H \in O(G)$ is non-monochromatic with $H \notin \mathcal{T} \cup \{D, \widetilde{D}\}$. Then both color classes in~$H$ form independent sets of size at least 2.
\end{lemma}

\begin{proof}
As before, we denote by $R_H$ and $B_H$ the red and the blue color class of $H$, respectively.
The possibilities for the structure of $R_H$ and $B_H$ are given by Theorem~\ref{theo:twoverticesnoedge}. Without loss of generality, we assume $|R_H| \geq |R_B|$. 
\medskip

We first consider the case that $R_H$ is a clique of size at least $3$ whose vertices are twins in $H$. By Theorem~\ref{theo:sizeclique}, we have $\omega_R= \aleph_0$.
If $B_H$ is a clique of size at least 3, then the vertices in $B_H$ are twins in $H$ by Theorem~\ref{theo:twoverticesnoedge} and, hence,~$H$ is homogeneously connected. Moreover, we have $\omega_B= \aleph_0$ by Theorem~\ref{theo:sizeclique}. By Lemma~\ref{lem:partition},~$H$ is realized in~$G[\rC \cup \bC]$ for any $\rC \in \mR$ and $\bC \in \mB$, which is a contradiction. Similarly, we argue if $H[B_H] \cong K_2$ and the vertices in $B_H$ are twins in $H$. Now assume $B_H = \{b_1, b_2\}$ is a 2-clique with $R_H \subseteq N_H(b_1)$ and $N_H(b_2) \cap R_H = \emptyset$. Since $|N_G^{\rC}(b)| = |\rC \setminus N_G^{\rC}(b)| = \aleph_0$ holds for all $b \in B$ and $\rC \in \mR$ (see Theorem~\ref{theo:sizeclique}), $H$ is realized in $G$, a contradiction. Now assume that $B_H$ is an independent set. Fix a red $|R_H|$-clique $K$ in $G$. By Lemma~\ref{lemma:clique+indepset}, we find an arbitrary number of pairwise non-adjacent joint blue neighbors and non-neighbors of $K$ in $G$. Hence $H$ is realized in $G$, which is a contradiction. 
\medskip

We proceed with the case that $H[R_H] \cong K_2$.
First assume that $|B_H| = 1$. Then $H \cong T_r$ or $H \cong \widetilde{T}_r$ (see Remark~\ref{rem:o3}), a contradiction.
Now assume $H[B_H] \cong K_2$. We obtain that $H$ is isomorphic to one of the graphs depicted in Figure~\ref{fig:proof}. 
\begin{figure}
	\centering
	
	\begin{subfigure}{1.2cm}
		\centering
				\begin{tikzpicture}[scale=.5, yscale=1.5]
					\node[\myred] (r1) at (0, 1) {};
					\node[\myred] (r2) at (0, 0) {};
					\node[\myblue] (b1) at (1.5, 1) {};
					\node[\myblue] (b2) at (1.5, 0) {};
					
					\draw[edge]
					(r1) edge (r2)
					edge (b1)
					edge (b2)
					(r2) edge (b1)
					(b1) edge (b2);
					\draw[edge] (r2) edge (b2);
				\end{tikzpicture}
	\end{subfigure}
	\hspace{.75cm}
	\begin{subfigure}{1.2cm}
		\centering
				\begin{tikzpicture}[scale=.5, yscale=1.5]
			\node[\myred] (r1) at (0, 1) {};
			\node[\myred] (r2) at (0, 0) {};
			\node[\myblue] (b1) at (1.5, 1) {};
			\node[\myblue] (b2) at (1.5, 0) {};
			
			\draw[edge]
			(r1) edge (r2)
			edge (b2)
			(r2) edge (b1) edge (b2);
			\draw[edge] (b1) edge (b2);
		\end{tikzpicture}
	\end{subfigure}
	\hspace{.75cm}
	\begin{subfigure}{1.2cm}
		\centering
			\begin{tikzpicture}[scale=.5, yscale=1.5]
			\node[\myred] (r1) at (0, 1) {};
			\node[\myred] (r2) at (0, 0) {};
			\node[\myblue] (b1) at (1.5, 1) {};
			\node[\myblue] (b2) at (1.5, 0) {};
			
			\draw[edge]
			(r2) edge (r1) edge (b1) edge (b2);
			\draw[edge]
			(b1) edge (b2);
	\end{tikzpicture}
	\end{subfigure}
	\hspace{.75cm}
	\begin{subfigure}{1.2cm}
		\centering
		\begin{tikzpicture}[scale=.5, yscale=1.5]
			\node[\myred] (r1) at (0, 1) {};
			\node[\myred] (r2) at (0, 0) {};
			\node[\myblue] (b1) at (1.5, 1) {};
			\node[\myblue] (b2) at (1.5, 0) {};
			
			\draw[edge]
			(r1) edge (r2) edge (b1);
			\draw[edge]
			(b2) edge (b1) edge (r2);
		\end{tikzpicture}
	\end{subfigure}
	\hspace{.75cm}
	\begin{subfigure}{1.2cm}
		\centering
			\begin{tikzpicture}[scale=.5, yscale=1.5]
				\node[\myred] (r1) at (0, 1) {};
				\node[\myred] (r2) at (0, 0) {};
				\node[\myblue] (b1) at (1.5, 1) {};
				\node[\myblue] (b2) at (1.5, 0) {};
				
				\draw[edge] (r1) to (r2) to (b2) to (b1);
			\end{tikzpicture}
	\end{subfigure}
	\hspace{.75cm}
	\begin{subfigure}{1.2cm}
		\centering
				\begin{tikzpicture}[scale=.5, yscale=1.5]
			\node[\myred] (r1) at (0, 1) {};
			\node[\myred] (r2) at (0, 0) {};
			\node[\myblue] (b1) at (1.5, 1) {};
			\node[\myblue] (b2) at (1.5, 0) {};
			
			\draw[edge]
			(r1) edge (r2);
			\draw[edge] (b1) edge (b2);
		\end{tikzpicture}	
	\end{subfigure}	
	\caption{Exhaustive list of possibilities of cross edges joining a red $K_2$ and a blue $K_2$.}\label{fig:proof}	
\end{figure}
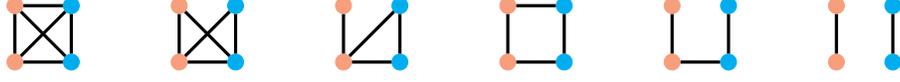
If $T_b$ and hence also $\widetilde{T}_b$ is omitted in~$G$ (see Lemma~\ref{lemma:nojointneighbors}), then $H$ is isomorphic to the fourth graph in Figure~\ref{fig:proof} since the other graphs are not minimally omitted.
In particular, all of the graphs in Figure~\ref{fig:proof} are omitted in $G$.
This is a contradiction since by the minimality of $H$, both the red $K_2$ and the blue $K_2$ are realized in~$G$ and, hence, at least one of the graphs given in Figure~\ref{fig:proof} is realized in $G$.
By symmetry, it follows that all graphs in $\mathcal{T}$ are realized in $G$. In particular, we have $\omega_R = \omega_B = \aleph_0$ (see Theorem~\ref{theo:sizeclique}).
Lemma~\ref{lem:partition} yields that all graphs in Figure~\ref{fig:proof} except for $D$ and $\widetilde{D}$ are realized in $G$. This is a contradiction. 
Finally, assume that $H[B_H] \cong \overline{K}_2$. %If $\omega_R= 2$, then~$H$ is realized in~$G$, a contradiction.
By Lemma~\ref{lem:partition} we have $\omega_R \neq 2$.
If $\omega_R = \aleph_0$, then $H$ is realized in $G$ by Lemma~\ref{lemma:clique+indepset}, which is a contradiction.

\medskip
The only remaining case is $|R_H| = |B_H| = 1$, which cannot occur since $G$ is not homogeneously connected.
\end{proof}

In the remaining part of this section, we show that the case described in Lemma~\ref{lemma:independentseth} does not occur. To this end, we need two technical lemmas. 

\begin{lemma}\label{lemma:technical}
Assume $\omega_R = \aleph_0$ and let $S \subseteq B$ be a finite independent set.
	For every $r \in R$, there exists $r' \in R$ with $r \sim r'$ and $N^S(r) = N^S(r')$.   
\end{lemma}

\begin{proof}
	Suppose that this is not the case. Let $H$ be the graph with vertex set $\{h_{x} \colon x \in S \cup \{r_1, r_2\}\}$ such that $h_{r_1} \sim h_{r_2}$ and for every $i,j \in \{1,2\}$ with $j \neq i$, the assignment $s \mapsto h_s$ for $s \in S$ and $r \mapsto r_i$ defines an isomorphism between $G[S \cup \{r\}]$ and $H-r_j$. Then $H$ is omitted in $G$. 
	First assume $\omega_B > 1$. If two adjacent red vertices in $G$ have a joint neighbor in a maximal blue clique, then they also have a joint non-neighbor (see Lemma~\ref{lem:partition}). Hence we may assume that all cross edges are present in $H$. But then $H$ is realized in $G$ by Lemma~\ref{lemma:clique+indepset}, which is a contradiction. Now assume $\omega_B = 1$ and let $n\coloneqq |N^S(r)|$. Omitting $H$ implies that for every $r'' \in R$ with $r'' \sim r$, the vertices~$r$ and $r''$ either have less than $n$ joint blue neighbors or less than $k-n$ joint blue non-neighbors which contradicts Lemma~\ref{lemma:clique+indepset}.
\end{proof}

\begin{lemma}\label{lemma:j}
Assume $\omega_R \geq 2$ and $\alpha_B = \aleph_0$. Let $J\subseteq R$ be a finite independent set. Fix $v \in B$ and set $S_v \coloneqq \{b \in B \colon N_G^J(b) = N_G^J(v)\}$. Then $S_v$ contains an infinite independent set.
\end{lemma}

\begin{proof}
We may assume $J \neq \emptyset$, so write $J = \{r_1, \ldots, r_t\}$ for some $t \in \N$ and $r_1, \ldots, r_t \in R$. 
For a contradiction, assume $\alpha(S_v) = k$ for some $k \in \N$. Then there exist $M_B^1, \ldots, M_B^k \in \mB$ with  $S_v \subseteq M_B^1 \cup \ldots \cup M_B^k \eqqcolon M$. Let $M_R \in \mR$ be the maximal red clique containing $r_t$.
	First assume $\omega_R = \aleph_0$. 
	For $r \in \rC$, let $S_r$ denote the set of vertices $s \in B$ for which the assignment $r_1 \mapsto r_1$, \dots, $r_t \mapsto r$, $v \mapsto s$ defines a partial isomorphism of~$G$. By ultrahomogeneity, we have $\alpha(S) = k$. By Lemma~\ref{lemma:technical}, we find $r' \in M_R$ with $N_G^{M}(r') = N_G^{M}(r)$. This implies $S = S_{r_t}\subseteq S_{r'}$ and hence $S = S_{r'}$ follows. Now assume that there exists $r'' \in \rC$ with $S_{r''} \not \subseteq M$.  
	Fixing $r_1, \ldots, r_t$ and exchanging $r'$ and $r''$ defines a partial isomorphism of~$G$, which permutes $\{M_B^1, \ldots, M_B^k\}$. On the other hand, it maps $S = S_{r'}$ to $S_{r''}$,
	and the latter set is not contained in $M$. This is a contradiction. Hence $S_{r''} \subseteq M$ holds for all $r'' \in M_R$.
	If $\omega_B= 1$ holds, this is impossible since every blue vertex has neighbors as well as non-neighbors in~$M_R$. For $\omega_B>1$, we use induction to find $M_B' \in \mB \setminus \{M_B^1, \ldots, M_B^k\}$ and a vertex $v' \in M_B'$ which has the same neighbors in $\{r_1, \ldots, r_{t-1}\}$ as~$v$. Since~$v'$ has both neighbors and non-neighbors in $M_R$, there exists $x \in M_R$ with $v' \in S_x$, which is a contradiction.
	\medskip
	
	Now assume $\omega_R = 2$. Let $r$ be the unique red neighbor of~$r_t$, and define the set $S_r$ as above. We argue by induction on $t$. Let $T_v$ be the set of blue vertices that have the same neighbors in $\{r_1, \ldots, r_{t-1}\}$ as $v$. By induction, $T_v$ contains an infinite independent set $U$. Then all but $k$ elements of $U$ lie in $S_r$. Hence $\alpha(S_r) = \aleph_0$ follows, which is a contradiction to ultrahomogeneity.
\end{proof}

With these results, we can now prove Theorem~\ref{theo:omittedgraphs}.

\begin{proof}[Proof of Theorem~\ref{theo:omittedgraphs}]
First assume that $D$ is omitted in $G$. Let $H \in O(G)$ be a non-monochromatic graph and assume $H \notin \mathcal{T} \cup \{D, \widetilde{D}\}$. By Lemma~\ref{lemma:independentseth}, the color classes in $H$ form independent sets of size at least~2. First assume $\omega_R \geq 2$ and $\alpha_B = \aleph_0$. Fix $b  \in B_H$. We view $H' \coloneqq H-b$ as induced subgraph of $G$. By minimality, $H[R_H \cup \{b\}]$ is realized in $G$, so the set $S \coloneqq \{s \in B \colon N_G^{R_H}(s) = N_G^{R_H}(b)\}$ is non-empty. By Lemma~\ref{lemma:j}, we find $s \in S$ non-adjacent to all blue vertices in $H'$. Then $G[V(H') \cup \{s\}]$ is isomorphic to $H$, a contradiction. If $\omega_R \geq 2$ holds and $\alpha_B$ is finite, we have $\omega_B = \alpha_R =\aleph_0$ by Theorem~\ref{theo:sizeclique}. Then we apply the above argument with interchanged colors. Similarly, we proceed if $\omega_R= 1$ holds since this implies $\alpha_R= \aleph_0$ and $\omega_B\geq 2$. 
%Now assume $\omega_R= 1$. Then $\alpha_R = \aleph_0$ and $\omega_B \geq 2$ follows. Hence we can apply Lemma~\ref{lemma:j} with exchanged colors and proceed as before. Similarly, we argue if $\alpha_B$ is finite since then $\alpha_R = \aleph_0$ follows by Theorem~\ref{theo:finitecliques}. 
Hence $H \in \mathcal{T} \cup \{D, \widetilde{D}\}$ follows. By Lemmas~\ref{lemma:nojointneighbors} and~\ref{lemma:h1h2}, we have $H \in O(G)$ if and only if $\widetilde{H} \in O(G)$ holds. Hence by Lemma~\ref{lemma:tildeomitted}, we obtain $G = \widetilde{G}$. 
\medskip

Now let $\widetilde{D}$ be omitted in $G$. By Lemma~\ref{lemma:tildeomitted}, the graph $\widetilde{G}$ omits $D$. The first part of this proof yields $G = \widetilde{\widetilde{G}} = \widetilde{G}$ and the claim follows.
\end{proof}

\subsection{\texorpdfstring{Classification of basic CUH graphs omitting $D$ or $\widetilde{D}$}{Classification of basic CUH graphs omitting D or Dtilde}}\label{sec:classification}
The first part of this section consists of examples of basic CUH graphs that omit $D$ and $\widetilde{D}$. The first one appears to be excluded by the unproven enumeration of possible CUH graphs stated in \cite{ROS11}.
Let $\mathcal{C}$ be the class of finite 2-colored graphs on red and blue vertices whose color classes form disjoint unions of cliques, that is, those graphs that omit monochromatic $P_3$s.

\begin{Example}\label{ex:1}
Consider the class $\mathcal{A}$ of all graphs in $\mathcal{C}$ which omit all of the following graphs:
a blue~$K_2$, a red~$K_3$, $T_r$, and $\widetilde{T_r}$.
% $B_A$ is an independent set and $R_A$ is a disjoint union of cliques of size at most $2$ that omit $T_r$ and $\widetilde{T}_r$.
We verify that $\mathcal{A}$ has the amalgamation property. Let $A_1, A_2 \in \mathcal{A}$, and let $J_1$ and $J_2$ be isomorphic induced subgraphs of $A_1$ and $A_2$, respectively. We construct an amalgam of $A_1$ and $A_2$ in $\mathcal{A}$ by identifying $J_2$ with the induced subgraph $J_1$ of $A_1$ and successively adding the vertices from $V(A_2) \setminus V(J_2)$ to~$A_1$. By induction, we may assume that $V(A_2) \setminus V(J_2)$ consists of a single vertex $v$. 

	\medskip
	
	First assume that $v$ is red and has a red neighbor $r \in V(J_1)$.
	If $r$ has a red neighbor $r'$ in~$A_1$ and $b \in V(J_1)$ is blue, then $r'\sim_{A_1}b$ if and only if $r \nsim_{J_1}b$, and this is the case precisely if $v\sim_{A_2} b$ since $T_r$ and $\widetilde{T}_r$ are omitted in $A_1$ and $A_2$. Hence mapping $V(J_1)$ to $V(J_2)$ and $r'$ to $v$ defines an isomorphism between $A_1[V(J_1) \cup \{r'\}]$ and $A_2$.
	If otherwise $N_{A_1}^R(r) = \emptyset$, then we obtain a graph $A$ from~$A_1$ by adding a new red vertex $r'$ with $N_A(r') = \{r\} \cup ( B_{A_1}\setminus N_{A_1}^B(r))$.
	% adjacent to $r$ and joining $r'$ to all blue non-neighbors of~$r$.
	Observe that $A$ omits~$T_r$ and $\widetilde{T}_r$ and contains $A_1$ as an induced subgraph.
	%In the following, we view $A_1$ % as an induced subgraph of $A$. 
	We claim that $A[(V(J_1) \cup \{r'\}]$ is isomorphic to $A_2$ in the canonical way. It suffices to check that for every blue vertex $b \in V(J_1)$, we have $b\sim_Ar'$ if and only if $b\sim_{A_2}v$.
	This is the case since $b\sim_{A_2}v$ precisely if $b\nsim_{A_2}r$, and this is the case if and only if $b\sim_Ar'$.
	
	\medskip
	Now assume that $v$ is red and $N_{J_1}^R(v) = \emptyset$.
	% does not have a red neighbor in~$J_1$.
	Consider the graph~$A$ obtained from $A_1$ by adding a new red vertex $r'$ with $N_A(r') = N_{J_1}^B(v)$.
	% which is adjacent to all blue vertices in $J_1$ to which $v$ is adjacent in $A_2$.
	It is easily verified that neither $T_r$ nor~$\widetilde{T}_r$ is realized in $A$ and $A[V(J_1) \cup \{r'\}] \cong A_2$.
	\medskip
	
	Finally assume that $v$ is blue. Again, we consider the graph $A$ obtained from $A_1$ by adding a blue vertex~$b$. We insert edges between $b$ and the red vertices in $A_1$ in two phases.
	First, we join $b$ to all vertices in $N^R_{A_2}(v) \cap V(J_1)$.
	% to all red vertices $w \in V(J_1)$ that are neighbors of $v$ in $A_2$.
	We mark the vertices in $R_{J_1} \setminus N_{A_2}(v)$.
	%All red vertices in $V(J_1)$ that are non-adjacent to $v$ in $A_2$ are marked.
	In the second phase, we iterate over all red cliques in $A_1$ that contain two vertices which are both not joined to $b$. For every such clique $\rC$, we add an edge joining $b$ with 
	an unmarked vertex in $\rC$. This is always possible: Assume that $\rC$ is a red clique in $A_1$ containing two marked vertices which are both non-neighbors of $b$ after the first phase. This implies $\rC \subseteq V(J_1)$, and both vertices in $\rC$ are non-adjacent to $v$ in $A_2$, which is a contradiction. In neither of the phases we generate an induced $T_r$, and 
	after the second phase, $A$ does not contain an induced $\widetilde{T}_r$.  
	Consider a red vertex $w \in V(J_1)$. If $w\sim_{A_2}v$, then~$wb$ is inserted in $A$ in the first phase. Otherwise, $w$ is marked in the first phase and hence $w\nsim_A b$. This shows that $A[V(J_1) \cup \{b\}]$ is isomorphic to~$A_2$. 
	\medskip
	
	Thus $\mathcal{A}$ is an amalgamation class. Its Fra\"issé limit $G$ satisfies $\omega_R(G) = 2$, and every blue vertex is joined to precisely one vertex of each maximal red clique. In particular, $G$ is not piecewise ultrahomogeneous.
\end{Example}

\begin{Example}\label{ex:2}
	Let $\mathcal{A}$ be the class of all graphs in $\mathcal{C}$ that omit monochromatic $K_3$s as well as the graphs in $\mathcal{T}$.
	 %$A$ in which $B_A$ and $R_A$ form disjoint unions of cliques of size at most two and that omit the graphs in $\mathcal{T}$. 
	 In order to show that $\mathcal{A}$ has the amalgamation property we use the approach and the terminology from Example~\ref{ex:1}. Let $A_1, A_2 \in \mathcal{A}$. By symmetry, we may assume that the vertex $v \in V(A_2) \setminus V(J_2)$ is red. 
	\medskip
	
	First assume that $v$ has a red neighbor $r$ in $A_2$. If $r$ has a red neighbor in $A_1$, we argue as in Example~\ref{ex:1} that $A_2$ is isomorphic to an induced subgraph of $A_1$. Otherwise, we add a red vertex $v'$ to $A_1$ which is adjacent to~$r$, and join $v'$ to all blue non-neighbors of $r$. Then every blue vertex has precisely one neighbor in $\{r, v'\}$. Assume that there exists a pair of adjacent blue vertices $b, b' \in V(A)$ which are both adjacent to $v'$. This means that they are both non-adjacent to $r$ in $A_1$, which is a contradiction. Now let $w \in V(J)$ be a blue vertex. We have $w\sim_{A_2}v$ if and only if $w\nsim_{A_2}r$, and this is equivalent to
	$w\sim_A v'$. Hence $A[V(J) \cup \{v'\}]$ is isomorphic to~$A_2$ in the canonical way. If $v$ does not have a red neighbor in~$A_2$, we add a red vertex $v'$ to $A_1$. We join $v'$ to the blue neighbors of $v$ in~$J$ and mark the blue non-neighbors of $v$ in $J$. For every blue $2$-clique $\bC$ in $A$ which does not contain neighbors of~$v'$, we then add an edge between $v'$ and an unmarked vertex in $\bC$. One can verify similarly to Example~\ref{ex:1} that this is possible, that the resulting graph is contained in $\mathcal{A}$, and that $A[V(J) \cup \{v'\}]$ is isomorphic to $A_2$. Hence~$\mathcal{A}$ is an amalgamation class.
\end{Example}

\begin{theorem}\label{theo:hiomittedclassification}
	Let $G$ be a basic CUH graph that omits $D$ or $\widetilde{D}$ and that is not isomorphic to one of the graphs in Theorem~\ref{theo:trussbipartite}. Up to interchanging the colors, one of the following cases arises:
	\begin{enumerate}[(i)]
		\item If $\alpha_R$ is finite, then $G$ is isomorphic to one of the following graphs: 
		\begin{enumerate}[(a)]
			\item The Fra\"issé limit $F_{\aleph_0,1}^{\alpha_R}$ of the class of graphs in $\mathcal{C}$ that omit the red $\overline{K}_{\alpha_R + 1}$ and the blue $K_2$.
			\item The Fra\"issé limit $F_{\aleph_0,2}^{\alpha_R}$ of the class of graphs in $\mathcal{C}$ that omit the red $\overline{K}_{\alpha_R + 1}$, the blue $K_3$ as well as $T_b$ and~$\widetilde{T}_b$. 
		\end{enumerate}
		\item Otherwise, $G$ is isomorphic to one of the following graphs: 
		\begin{enumerate}[(a)]
			\item The Fra\"issé limit $F_{2,1}$ of the class of graphs in $\mathcal{C}$ that omit the red $K_3$, the blue $K_2$ as well as $T_r$ and~$\widetilde{T}_r$. 
			\item The Fra\"issé limit $F_{\aleph_0,1}$ of the class of graphs in $\mathcal{C}$ that omit the blue $K_2$.
			\item The Fra\"issé limit $F_{2,2}$ of the class of graphs in $\mathcal{C}$ that omit monochromatic $K_3$s and the graphs in~$\mathcal{T}$.
			\item The Fra\"issé limit $F_{\aleph_0,2}$ of the class of graphs in $\mathcal{C}$ that omit the blue $K_3$, $T_b$, and $\widetilde{T}_b$.
			\item The Fra\"issé limit $F_{\aleph_0, \aleph_0}$ of the class of graphs in $\mathcal{C}$ that omit $D$ and $\widetilde{D}$.
		\end{enumerate}
	\end{enumerate}
\end{theorem}
\begin{proof}
First, we verify that the given classes are amalgamation classes. For the cases (a) and (c) in the second part this is done in Examples~\ref{ex:1} and~\ref{ex:2}, respectively. The remaining cases can be treated similarly. Now let~$G$ be a basic CUH graph omitting $D$ or $\widetilde{D}$ that is not isomorphic to one of the graphs in Theorem~\ref{theo:trussbipartite}.
By Theorem~\ref{theo:sizeclique}, we have $\omega_R,\omega_B \in \{1,2,\aleph_0\}$. Without loss of generality, we assume $\omega_R \geq \omega_B$. The case $\omega_R = \omega_B = 1$ is excluded since $G$ is not isomorphic to one of the graphs in Theorem~\ref{theo:trussbipartite}. If $\omega_R = 2$ and $\omega_B = 1$, then $G$ minimally omits $T_{r}$ and $\widetilde{T}_{r}$ by Theorem~\ref{theo:omittedgraphs} and Lemma~\ref{lemma:nojointneighbors}, which yields $G \cong F_{2,1}$. For $\omega_R = \omega_B = 2$, the graphs in $\mathcal{T}$ are minimally omitted in $G$ and, hence, we have $G \cong F_{2,2}$. If $\omega_R = \aleph_0$ and $\omega_B = 1$, we obtain $G \cong F_{\aleph_0,1}$ or $G \cong F_{\aleph_0,1}^{\alpha_R}$, depending on whether $\alpha_R$ is infinite or finite. 
%Note that the case $\alpha_R = 1$ is excluded by the assumption that $G$ is basic. 
If $\omega_R = \aleph_0$ and $\omega_B = 2$, then $G$ omits $T_{b}$ and $\widetilde{T}_{b}$ by Lemma~\ref{lemma:nojointneighbors}. Depending on whether $\alpha_R$ is infinite or finite, we 
obtain $G \cong F_{\aleph_0,2}$ or $G \cong F_{\aleph_0,2}^{\alpha_R}$. Finally, if $\omega_R = \omega_B = \aleph_0$ holds, we have $\alpha_R = \alpha_B = \aleph_0$ by Theorem~\ref{theo:finitecliques}. Moreover, all graphs in $\mathcal{T}$ are realized in~$G$ by Lemma~\ref{lemma:nojointneighbors}. By Lemma~\ref{lemma:h1h2}, we have $O(G) = \{D, \widetilde{D}\}$ and, hence, $G \cong F_{\aleph_0, \aleph_0}$.	
\end{proof}

Combining Theorems~\ref{theo:bedingungenuh} and~\ref{theo:hiomittedclassification} yields the characterization of piecewise ultrahomogeneity stated in Theorem~\ref{theo:a}.

\section{Classification of piecewise ultrahomogeneous graphs}\label{sec:cliquereducible}
In this section, we complete the proof of Theorem~\ref{theo:b} by classifying the basic piecewise ultrahomogeneous CUH graphs. Unless stated otherwise, we assume that $G$ is a basic piecewise ultrahomogeneous CUH graph that is neither isomorphic to $F_{2,2}$ nor to one of the graphs in Theorem~\ref{theo:trussbipartite}. The last condition is equivalent to requiring that $D$ and~$\widetilde{D}$ are realized in $G$ (see Theorem~\ref{theo:a}).

\begin{Remark}\label{rem:allequal}
Since $G$ is not a blow-up, there exist both edges and non-edges between every $\rC \in \mR$ and $\bC \in \mB$. This yields $G[\rC \cup \bC] \cong G[\rC' \cup \bC']$ for all $\rC,\rC' \in \mR$ and $\bC,\bC' \in \mB$.
\end{Remark}

\begin{lemma}\label{lemma:perfectmatching}
For every $\rC \in \mR$ and $\bC \in \mB$, the graph $G[\rC \cup \bC]$ is isomorphic to the complement of the generic bipartite graph. 
\end{lemma}

\begin{proof}
The graph $\overline{G[\rC \cup \bC]}$ is an ultrahomogeneous graph whose color classes form independent sets. By Theorem~\ref{theo:trussbipartite}, $\overline{G[\rC \cup \bC]}$ is isomorphic to the generic bipartite graph since $D$ and $\widetilde{D}$ are realized in $G$ and, hence, in $G[\rC \cup \bC]$ by Remark~\ref{rem:allequal}.
\end{proof}

In particular, we have $\omega_R = \omega_B= \aleph_0$. We now prove that the isomorphism type of $G$ only depends on $\alpha_R$ and $\alpha_B$. To this end, we need the following technical lemma:

\begin{lemma} \label{lemma:jointneighborgeneric}
For every $U \subseteq R$ and every $\bC \in \mB$ there exists a vertex $b \in \bC$ with $U \subseteq N(b)$.
%Let $K \subseteq V(G)$ be a finite set of red vertices. Every maximal blue clique $B \in \mB$ contains a joint neighbor of $K$. 
\end{lemma}

\begin{proof}
We fix $U \subseteq R$ and $\bC \in \mB$.
We have $N_G^B(U) \neq \emptyset$ since the graph arising from~$G[U]$ by adding a blue vertex dominating $U$ cannot be omitted in $G$ (using Theorem~\ref{theo:twoverticesnoedge} and Lemma~\ref{lemma:perfectmatching}). Suppose that no vertex in $\bC$ dominates $U$. We fix a vertex $v \in \bC$. Let $H$ be a finite graph with $V(H) = \{h_x \colon x \in U \cup \{v,v'\}\}$ such that the assignment $h_x \mapsto x$ defines an isomorphism $\varphi$ between $H-h_{v'}$ and $G[U \cup \{v\}]$, and such that the blue vertex $h_{v'}$ is adjacent to all vertices in $H$. Assume that $H$ is realized in $G$ and let $\hat{\varphi} \in \Aut(G)$ be an extension of $\varphi$ to $G$. Then $\hat{\varphi}(v')$ dominates $U$ and is contained in $\bC$. This is a contradiction, so $H$ is omitted in $G$. Let $H'$ be an induced subgraph of $H$ that is minimally omitted in $G$. Then $H'$ contains both $h_v$ and $h_{v'}$ since $H-h_v$ and $H-h_{v'}$ are realized in~$G$ by assumption. By Theorem~\ref{theo:twoverticesnoedge}, $R_{H'}$ is an independent set of size at least 2 since $H'$ is realized in $G[\rC \cup \bC]$ for any $\rC \in  \mR$ otherwise. Fix two adjacent blue vertices $b, b' \in V(G)$. By Lemma~\ref{lemma:perfectmatching}, every maximal red clique of $G$ contains vertices which are adjacent to both, to precisely one, or to none of the vertices in $\{b,b'\}$. Hence $H'$ cannot be omitted in $G$, which is a contradiction.
\end{proof}

Of course, the analogous statement with interchanged colors holds as well.
\medskip

	For every $r, b \in \N \cup \{\aleph_0\}$, let $\mathcal{A}_{r,b}$ denote the class of finite 2-colored graphs in which the red vertices form disjoint union of at most $r$ cliques, and the blue vertices form a disjoint union of at most $b$ cliques. It is easily verified that $\mathcal{A}_{r,b}$ is an amalgamation class. Let $G_{r,b}$ denote its Fra\"issé limit.

\begin{lemma}\label{lemma:gb}
The graph $G$ is isomorphic to $G_{\alpha_R,\alpha_B}$.
\end{lemma}

\begin{proof}
Let $\mathcal{A}_G$ be the age of $G$. Clearly, we have $\mathcal{A}_G \subseteq \mathcal{A}_{\alpha_R,\alpha_B}$. Now let $H \in \mathcal{A}_{\alpha_R,\alpha_B}$ and suppose for a contradiction that $H$ is omitted in $G$. We may assume $H \in O(G)$. It is easy to see that $H$ is not monochromatic. By Lemma~\ref{lemma:perfectmatching}, one of the color classes in $H$ does not form a clique. By Theorem~\ref{theo:twoverticesnoedge}, we may assume that the blue vertices in $H$ form an independent set of size at least two. Fix a blue vertex $b \in V(H)$. By assumption, $H' \coloneqq H-b$ is realized in~$G$, and we view it as an induced subgraph of $G$. Let $K$ be an induced subgraph of~$G$ on vertices $V(K) = \{k_x \colon x \in R_{H'} \cup \{b\}\}$ such that the assignment $k_x \mapsto x$ defines an isomorphism from $K$ to $H[R_{H'} \cup \{b\}]$. Let $\bC' \in \mB$ be the maximal blue clique of $G$ containing $k_b$. There exists $\bC \in \mB$ not containing a vertex of $H'$, as $H'$ contains at most $\alpha_B-1$ blue vertices. By Lemma~\ref{lemma:jointneighborgeneric}, there exist vertices $v \in \bC$ and $v' \in \bC'$ dominating $R_{H'} \cup R_K$. Mapping~$R_K$ bijectively to $R_{H'}$ and $v'$ to $v$ defines a partial isomorphism $\varphi$ of $G$. Let $\hat{\varphi} \in \Aut(G)$ be an extension of $\varphi$ to $G$. Then $\hat{\varphi}(k_b) \in \bC$ has the same neighbors in $R_{H'}$ as $b$. Due to $\bC \cap V(H') = \emptyset$, extending $H'$ by $\hat{\varphi}(k_b)$ gives rise to an embedding of $H$ into $G$. This is a contradiction. Hence we have $\mathcal{A}_G = \mathcal{A}_{\alpha_R, \alpha_B}$. By Theorem~\ref{theo:fraisse}, the graphs $G$ and $G_{\alpha_R, \alpha_B}$ are isomorphic. 
\end{proof}

\begin{Remark}
The graph $G_{\alpha_R, \alpha_B}$ is generic in the following sense: For every color $c \in \{\text{red, blue}\}$, every maximal clique $C$ of color~$c$ and all finite disjoint vertex sets $S$ and $T$ of color $c' \neq c$, there exists a vertex $v \in C$ with $S \subseteq N_G(v)$ and $N_G(v) \cap T = \emptyset$.
\end{Remark}

Summarizing, we obtain the following classification of basic piecewise ultrahomogeneous CUH graphs:

	\begin{theorem}\label{theo:twocolors}
		Let $G$ be a basic piecewise ultrahomogeneous CUH graph. Then $G$ is isomorphic to one of the graphs in Theorem~\ref{theo:trussbipartite}, to $F_{2,2}$, or to~$G_{\alpha_R,\alpha_B}$. 
	\end{theorem}

This completes the proof of Theorem~\ref{theo:b}. 

\section{Conclusion}\label{sec:conclusion}

In this paper, we classified the countable ultrahomogeneous 2-colored graphs in which each color class forms a disjoint union of cliques.
Since complementing a color class preserves ultrahomogeneity, this directly translates to a full classification of ultrahomogeneous 2-colored graphs with imprimitive color classes.
 Our key tool was the concept of piecewise ultrahomogeneity introduced in Section~\ref{sec:neighboringrelations}. We showed that with one exception, a basic non-bipartite CUH graph is piecewise ultrahomogeneous if and only if two specific graphs appear as induced subgraphs (see Theorem~\ref{theo:a}). Using this result, we obtained the classification of countable 2-colored CUH graphs given in Theorem~\ref{theo:b}.
The existence of non-piecewise ultrahomogeneous graphs is a strong contrast to the finite case, where every basic ultrahomogeneous graph is also piecewise ultrahomogeneous~\cite{HEI20}.
\medskip

There are several natural continuations of this paper. For example, it would be interesting to classify edge-colored versions of CUH graphs, extending the work of Lockett and Truss \cite{LOC14}. Moreover, one could investigate $n$-colored versions of CUH graphs for an arbitrary number $n \in \N$. In both cases, we believe that a suitable generalization of piecewise ultrahomogeneity could play a central role. Just as in the case studied in this paper, one could hope to characterize the piecewise ultrahomogeneous graphs in terms of a small number of induced subgraphs, and then use the classifications of ultrahomogeneous multipartite graphs given in \cite{JEN12} and \cite{LOC14}. Conversely, if a graph fails to be piecewise ultrahomogeneous, its structure might again be very limited.

\bibliographystyle{plain}
\bibliography{cliquescombined.bib}
\end{document}